\documentclass[11pt,bezier]{article}
\usepackage{amsmath,amssymb,amsfonts,euscript}

\usepackage{algorithmic, algorithm}

\usepackage{graphicx}

\usepackage[ansinew]{inputenc}
\usepackage{graphicx}
\usepackage{color}
\usepackage{mathrsfs}

\usepackage[colorlinks]{hyperref}
\usepackage[active,new,noold,marker]{xrcs}
\catcode`\=13
\def{$\bowtie$}

\usepackage{eurosym}
\usepackage{lscape} 

\newtheorem{prethm}{{\bf Theorem}}

\newenvironment{thm}{\begin{prethm}{\hspace{-0.5
               em}{\bf}}}{\end{prethm}}

\newtheorem{prepro}{{\bf Theorem}}

\newenvironment{pro}{\begin{prepro}{\hspace{-0.5
               em}{\bf}}}{\end{prepro}}

\newtheorem{preprop}{{\bf Proposition}}

\newtheorem{precor}{{\bf Corollary}}

\newtheorem{preconj}{{\bf Conjecture}}

\newtheorem{predefi}{{\bf Definition}}

\newtheorem{preremark}{{\bf Remark}}

\newtheorem{preexample}{{\bf Fact}}

\newtheorem{prelem}{{\bf Lemma}}

\newtheorem{prelam}{{\bf Lemma}}

\newtheorem{preprob}{{\bf Problem}}

\newenvironment{prob}{\begin{preprob}{\hspace{-0.5
               em}{\bf.}}}{\end{preprob}}

\newtheorem{preproof}{{\bf Proof}}

\newtheorem{preali}{{\bf Proof of Theorem 1.}}

\newenvironment{ali}[1]{\begin{preali}{\rm
               #1}\hfill{$\Box$}}{\end{preali}}

\newtheorem{prealii}{{\bf Proof of Theorem 2.}}

\newenvironment{alii}[1]{\begin{prealii}{\rm
               #1}\hfill{$\Box$}}{\end{prealii}}

\newtheorem{prealiii}{{\bf Proof of Theorem 3.}}

\newenvironment{aliii}[1]{\begin{prealiii}{\rm
               #1}\hfill{$\Box$}}{\end{prealiii}}

\newtheorem{prealiiii}{{\bf Proof of Theorem 4.}}

\newenvironment{aliiii}[1]{\begin{prealiiii}{\rm
               #1}\hfill{$\Box$}}{\end{prealiiii}}

\newtheorem{prealiiiii}{{\bf Proof of Theorem 5.}}

\newenvironment{aliiiii}[1]{\begin{prealiiiii}{\rm
               #1}\hfill{$\Box$}}{\end{prealiiiii}}

\newtheorem{prealiiiiii}{{\bf Proof of Theorem 6.}}

\newenvironment{aliiiiii}[1]{\begin{prealiiiiii}{\rm
               #1}\hfill{$\Box$}}{\end{prealiiiiii}}

\newtheorem{prealiiiiiii}{{\bf Proof of Theorem 7.}}

\newenvironment{aliiiiiii}[1]{\begin{prealiiiiiii}{\rm
               #1}\hfill{$\Box$}}{\end{prealiiiiiii}}

\newtheorem{prealij}{{\bf Proof of Theorem 8.}}

\newenvironment{alij}[1]{\begin{prealij}{\rm
               #1}\hfill{$\Box$}}{\end{prealij}}

\newtheorem{prealijj}{{\bf Proof of Theorem 9.}}

\newenvironment{alijj}[1]{\begin{prealijj}{\rm
               #1}\hfill{$\Box$}}{\end{prealijj}}

\title{On the algorithmic complexity of adjacent vertex closed distinguishing colorings
number of graphs}

\author{{\normalsize
{  Ali Dehghan${}^{\mathsf{a}}$},\,
{  Mohsen Mollahajiaghaei${}^{\mathsf{b}}$},\,
}\vspace{3mm}
\\{\footnotesize{${}^{\mathsf{a}}$\it Systems and Computer Engineering Department, Carleton University, Ottawa,   Canada}}  {\footnotesize{}}\\{\footnotesize{${}^{\mathsf{b}}$\it
Department of Mathematics, University of Western Ontario, London, Ontario, Canada}}
\thanks{{\it E-mail addresses}:  $\mathsf{alidehghan@sce.carleton.ca}$, $\mathsf{mmollaha@uwo.ca}$.} }

\date{}
\begin{document}
\maketitle

\begin{abstract}
{\small \noindent
An assignment of numbers  to the
vertices of graph $G$ is {\it closed distinguishing} if for any two adjacent vertices $v$ and $u$ the sum of labels of the vertices in the closed neighborhood of the vertex $v$ differs from
the sum of labels of the vertices in the closed neighborhood of the vertex $u$  unless they have the same closed neighborhood (i.e. $N[u] = N[v]$).
The {\it closed distinguishing number} of a graph $G$, denoted by $dis[G]$, is
the smallest integer $k$ such that there is a closed distinguishing labeling for
$G$ using integers from the set $\{1,2,\ldots,k\}$.
Also, for each vertex $v \in V(G)$, let $L(v)$ denote a list of natural numbers available at $v$. A  {\it list closed distinguishing labeling} is a closed distinguishing labeling $f$ such that $f(v) \in L(v)$ for each $v \in V(G)$.
A graph $G$ is said to be {\it closed distinguishing $k$-choosable} if every $k$-list assignment of natural numbers to the vertices of $G$ permits a list closed distinguishing labeling of $G$. The closed distinguishing choice number of $G$, $dis_{\ell}[G]$, is the minimum natural number $k$ such that $G$ is closed distinguishing $k$-choosable.
In this work we show that for each integer $t$ there is a bipartite graph $G$ such that $dis[G] > t$. This is an answer to a question raised by Axenovich et al. in \cite{dis} that how ”dis” function depends on the chromatic number of a graph.
It was shown that for every graph $G$ with $\Delta\geq 2$, $dis[G]\leq dis_{\ell}[G]\leq \Delta^2-\Delta+1$ and also there are infinitely many values of $\Delta$ for which $G$ might be chosen so that $dis[G] =
\Delta^2-\Delta+1$ \cite{dis}. In this work, we prove that the difference between $dis[G]$ and $dis_{\ell}[G]$ can be arbitrary large and show that for every positive integer $t$ there is a graph $G$ such that $dis_{\ell}[G]-dis[G]\geq t$.
Also, we improve the current upper bound and give some number of upper bounds for the closed distinguishing choice number by using the Combinatorial Nullstellensatz.
Among other results, we show that
it is $ \mathbf{NP} $-complete to decide for a given planar subcubic graph $G$, whether $ dis[G]=2$.
Also, we prove that for every $k\geq 3$, it is
{\bf NP}-complete to decide whether $dis[G]= k$ for a given graph
$G$.
}

\begin{flushleft}


\end{flushleft}

\end{abstract}

\section{Introduction}
\label{}
In 2004, Karo\'nski et al. in \cite{MR2047539} introduced a new coloring of
a graph which is generated via edge labeling. Let $f : E(G) \rightarrow \mathbb{N}$ be a labeling of the edges of a graph $G$ by positive integers such that for every two adjacent vertices $v$ and $u$, $S(v) \neq S(u)$, where $S(v)$ denotes the sum of labels of all edges
incident with $v$. It was conjectured that three integer labels
$\{1,2,3\}$ are sufficient for every connected graph, except
$K_2$ \cite{MR2047539} (1-2-3 Conjecture). Currently the best bound that was proved by Kalkowski et al. is five \cite{MR2595676}. For more information we refer the reader to  a survey on the 1-2-3 Conjecture and related problems by Seamone \cite{survey} (also see \cite{bartnicki2009weight, MR3512668, MR3315373, MR3072733, david,  MR2654258,  wong2012total}). Different
variations of distinguishing labelings of graphs have also been considered, see \cite{baudon2015oriented, bensmail2014partitions, gyHori2009new, khatirinejad2011digraphs, khatirinejad2012vertex,  lu2011vertex, seamone2012derived, seamone2014bounding,  seamone2012sequence}.

On the other hand, there are different types of labelings which consider the closed neighborhoods of vertices.
In 2010, Esperet et al. in  \cite{esperet2010locally} introduced the notion of locally identifying coloring of a graph.
A proper vertex-coloring of a graph $G$ is said to be {\it locally identifying} if for any pair $u$, $v$ of adjacent vertices with distinct closed neighborhoods, the sets of colors in the closed neighborhoods of $u$ and $v$ are different.
In 2014, A{\"\i}der et al. in  \cite{aider2014relaxed} introduced the notion of relaxed locally identifying coloring of graphs.
A  vertex-coloring of a graph $G$ (not necessary proper) is said to be {\it relaxed locally identifying} if for any pair $u$, $v$ of adjacent vertices with distinct closed neighborhoods, the sets of colors in the closed neighborhoods of $u$ and $v$ are different.
Note that a relaxed locally identifying coloring of
a graph that is similar to locally identifying coloring for which the coloring
is not necessary proper.
For more information see \cite{ foucaud2012locally, gonccalves2013locally, survey}.

Motivated by the 1-2-3 Conjecture and the relaxed locally identifying coloring, the closed distinguishing labeling as a vertex version of the 1-2-3 Conjecture was introduced  by Axenovich et al. \cite{dis}.
For every vertex $v$ of $G$, let $N[v]$ denote the closed neighborhood of
$v$.
An assignment of numbers to the
vertices of a graph $G$ is {\it closed distinguishing}
if for any two adjacent vertices $v$ and $u$ the sum of labels of the vertices in the closed neighborhood of the vertex $v$ differs from
the sum of labels of the vertices in the closed neighborhood of the vertex $u$
unless  $N[u] = N[v]$ (i.e. they have the same closed neighborhood).
The {\it closed distinguishing number} of a graph $G$, denoted by $dis[G]$, is
the smallest integer $k$ such that there is a closed distinguishing assignment for
$G$ using integers from the set $\{1,2,\ldots,k\}$.
For each vertex $v \in V(G)$, let $L(v)$ denote a list of natural numbers available at $v$. A  {\it list closed distinguishing labeling} is a closed distinguishing labeling $f$ such that $f(v) \in L(v)$ for each $v \in V(G)$.
A graph $G$ is said to be {\it closed distinguishing $k$-choosable} if every $k$-list assignment of natural numbers to the vertices of $G$ permits a list closed distinguishing labeling of $G$. The {\it closed distinguishing choice number} of $G$, $dis_{\ell}[G]$, is the minimum natural number $k$ such that $G$ is closed distinguishing $k$-choosable. In this work we study
closed distinguishing number and closed distinguishing choice number of graphs.

The closed distinguishing number of a graph $G$ is
the smallest integer $k$ such that there is a closed distinguishing assignment for
$G$ using integers from the set $\{1,2,\ldots,k\}$. In this work, we also consider another parameter, the minimum number of integers
required in a closed distinguishing labeling. For a given graph $G$, the minimum number of integers
required in a closed distinguishing labeling is called its {\it strong closed distinguishing number} $dis_s[G]$.
Note that a  vertex-coloring of a graph $G$ (not necessary proper) is said to be strong closed distinguishing labeling  if for any pair $u$, $v$ of adjacent vertices with distinct closed neighborhoods, the multisets of colors in the closed neighborhoods of $u$ and $v$ are different.

\section{Closed distinguishing labeling}

In this section we study closed distinguishing number and closed distinguishing choice number of graphs.

\subsection{The difference between $dis[G]$ and $dis_{\ell}[G]$}

It was shown in  \cite{dis} that for every graph $G$ with $\Delta\geq 2$, $dis[G]\leq dis_{\ell}[G]\leq \Delta^2-\Delta+1$.
Also, there are infinitely many values of $\Delta$ for which $G$ might be chosen so that $dis[G] =
 \Delta^2-\Delta+1$ \cite{dis}.
We prove that the difference between $dis[G]$ and $dis_{\ell}[G]$ can be arbitrary large and show that for every number $t$ there is a graph $G$ such that $dis_{\ell}[G]-dis[G]\geq t$.

\begin{thm}\label{T1}
For every positive integer $t$ there is a graph $G$ such that $dis_{\ell}[G]-dis[G]\geq t$.
\end{thm}

\subsection{The complexity of determining $dis[G]$}

Let $T \neq K_2$ be a tree. It was shown  \cite{dis} that $dis_{\ell}[T]\leq 3$ and $dis[T]\leq 2$.
Here, we investigate the computational complexity of determining   $dis[G]$ for planar subcubic graphs and bipartite subcubic graphs.

\begin{thm}\label{T2}
For a given planar subcubic graph $G$, it is $ \mathbf{NP} $-complete to decide  whether $ dis[G]=2$.
\end{thm}

Although for a given tree $T$, we can compute $dis[T]$ in polynomial time \cite{dis}, but the problem  of determining the closed distinguishing number is hard for bipartite graphs.

\begin{thm}\label{T3}
 For a given bipartite subcubic graph $G$, it is $ \mathbf{NP} $-complete to decide  whether $ dis[G]=2$.
\end{thm}

Note that in the proof of Theorem \ref{T3}, we reduced Not-All-Equal to our problem and the planar version of Not-All-Equal
is in $ \mathbf{P} $ \cite{p}, so the computational complexity of
deciding whether $ dis[G]=2$ for planar bipartite graphs remains unsolved.

\begin{thm}\label{NT2}
For every integer $t\geq 3$, it is
{\bf NP}-complete to decide whether $dis[G]= t$ for a given graph
$G$.
\end{thm}

\subsection{Upper bounds for $dis_{\ell}[G]$ and $dis[G]$}

It was shown that for every graph $G$ with $\Delta\geq 2$, $dis[G]\leq dis_{\ell}[G]\leq \Delta^2-\Delta+1$ \cite{dis}.
Here, we improve the previous bound.

\begin{thm}\label{T4}
Let $G$ be a simple graph on $n$ vertices with degree sequence $\Delta=d_{1}\geq d_{2}\geq \cdots \geq d_{n}=\delta$ and $\Delta\neq 1$. Define $s:=d_1+\cdots+d_{\Delta}-\Delta$.\\
(i) $\displaystyle   dis_{\ell}[G] \leq s+1 \leq \Delta^2-\Delta+1$.\\
(ii) $dis_{\ell}[G]\leq m$, where $m$ is the number of edges.\\
(iii) If there are exactly $t$ vertices with  degree $\Delta$. Then
\begin{center}
$\displaystyle    dis_{\ell}[G] \leq \min\{\Delta^{2}-2\Delta +t +1,\Delta^2-\Delta+1\}$.
\end{center}
(iv) If there is a unique vertex with  degree $\Delta$. Then
$\displaystyle   dis_{\ell}[G] \leq \Delta^2-3\Delta+4$. \\
(v) If $G$ is a strongly regular graph with parameters $(n,k,\lambda,\mu)$. Then
\begin{center}
$\displaystyle  dis_{\ell}[G] \leq k(k-\lambda-1) +1$.
\end{center}
(vi) $\displaystyle   dis_{\ell}[G] \leq (\dfrac{n-1}{2})^{2}+1$.
\end{thm}

\subsection{Lower bound for $dis[G]$}

Let $G$ be a bipartite graph with partite sets $A$ and $B$ which is not a star. Let,
for $X\in \{A,B \}$; $\displaystyle \Delta_{X} = \max_{x\in X}d(x) $
and $\displaystyle\delta_{X,2} = \min_{x\in X, d(x)\geq 2}d(x)$. It was shown in \cite{dis} that
$$\displaystyle dis[G] \leq min\{c\sqrt{|E(G)|},\lfloor \dfrac{\Delta_A -1}{\delta_{B,2}-1}\rfloor+1,  \lfloor \dfrac{\Delta_B -1}{\delta_{A,2}-1}\rfloor+1\}, $$
where $c$ is some constant. Thus, for a given bipartite graph $G$, $dis[G]=\mathcal{O}(\Delta)$ \cite{dis}.
Regarding "dis" function, Axenovich et al. in \cite{dis} said: "One of the challenging problems in the area is to determine how "dis" function depends
on the chromatic number of a graph. The situation is far from being understood even for
bipartite graphs."
We give a negative answer to this problem and show that for each $t$ there is a bipartite graph $G$ such that $dis[G]>t$.

\begin{thm}\label{Tj1}
For each integer $t$, there is a bipartite graph $G$ such that $dis[G]>t$.
\end{thm}

\subsection{Split graphs}

It is well-known that split graphs can be recognized in polynomial time, and that finding a canonical partition of a split graph can also be found in polynomial time. We prove the following result.

\begin{thm}
If $G$ is a split graph, then $dis[G]\leq(\omega(G))^2$.
\end{thm}

\section{Strong closed distinguishing number}

In this section, we focus on the strong closed distinguishing number of graphs.
For any graph $G$, we have the following.

\begin{equation}
 dis_s[G] \leq dis[G] \leq dis_{\ell}[G]
\end{equation}
For a given bipartite graph $G=[X,Y]$, define $f:V(G)\rightarrow \{1,\Delta\}$ such that:\\
$$ f(v)=
\begin{cases}
   1,           &$     if  $ v\in X \\
   \Delta,    &$     if  $ v\in Y
\end{cases}$$
It is easy to see that $f$ is a closed distinguishing labeling for $G$. Thus, $dis_s[G]\leq 2 $.
So, by Theorem \ref{Tj1}, the difference between $dis[G]$ and $dis_{s}[G]$ can be arbitrary large.
Here we increase the gap.

\begin{thm}\label{T6}
For each $n$, there is a graph  $G$ with $n$ vertices such that $dis[G]-dis_{s}[G]=\Omega(n^{\frac{1}{3}})$.
\end{thm}

Let $G$ be an $r$-regular graph and $f:V(G)\rightarrow \{a,b\}$ be a closed distinguishing labeling.
Define:
$$g(v)=
\begin{cases}
   a',      &$     if  $\,\,  f(v)=a, \\
   b',      &$     if  $\,\,  f(v)=b.
\end{cases}
$$
It is easy to check that if $a'\neq b'$, then $g:V(G)\rightarrow \{a',b'\}$ is a closed distinguishing labeling.
Thus, for an $r$-regular graph $G$,  $dis_{s}[G]=2$ if and only if $dis[G]=2$.

Let $a$ and $b$ be two numbers and $a\neq b$,  we show that for a given 4-regular graph $G$, it is $ \mathbf{NP} $-complete to decide  whether there is a closed distinguishing labeling from $\{a,b\}$.

\begin{thm}\label{T7}
For a given 4-regular graph $G$, it is $ \mathbf{NP} $-complete to decide  whether $ dis_{s}[G]=2$.
\end{thm}

\section{Notation and Tools}

All graphs considered in this paper are finite, undirected, with no loops or multiple
edges. If $G$ is a graph, then $V(G)$ and $E(G)$ denote the vertex set and the edge set of
$G$, respectively. Also, $\Delta(G)$ denotes the maximum degree of $G$ and simply denoted
by $\Delta$. For every $v\in V(G)$, $d_{G}(v)$ and $N_{G}(v)$ denote the degree of $v$ and the set of
neighbors of $v$, respectively. Also $N[v]=N(v)\cup \{v\}$.
For a given graph $G$, we use $u\sim v$ if two vertices $u$ and $v$ are adjacent in
$G$.

Let $G$ be a graph and $K$ be a non-empty set. A proper vertex coloring of $G$ is a
function $c : V(G)\rightarrow K$, such that if $u, v \in V(G)$ are adjacent, then $c(u)\neq c(v)$.
A {\it proper vertex $k$-coloring} is a proper vertex coloring with $|K| = k$. The smallest number
of colors needed to color the vertices of $G$ for obtaining a proper vertex coloring is
called the {\it chromatic number} of $G$ and denoted by $\chi(G)$.

A {\it $k$-regular} graph is a graph whose each vertex has degree $k$.
A regular graph $G$ with $n$ vertices and degree $k$ is said to be {\it strongly regular} if there
are integers $ \lambda $ and $ \mu $ such that every two adjacent vertices have $ \lambda $ common neighbors
and every two non-adjacent vertices have $ \mu $ common neighbors and is denoted by
$SRG(n,k,\lambda,\mu)$.

The {\it
Cartesian product} $H\square G$ of graphs $G$ and $H$ is the graph with vertex set $V(G)\times V(H)$
where vertices $(g, h)$ and $(g', h')$ are adjacent if and only if either $g = g'$ and $hh'\in E(H)$, or $h = h'$ and $gg'\in E(G)$.

We say that a set of vertices is {\it independent} if there is no edge between
these vertices. The {\it independence number}, $\alpha(G)$, of a graph $G$ is the size of a largest independent set of $G$.
Also, a {\it clique} in a  graph $G$  is a subset of its vertices such that every two vertices in the subset are connected by an edge. The {\it clique number} $\omega(G)$ of a graph $G$ is the number of vertices in a maximum clique in $G$.
A {\it split graph} is a graph whose vertex set may be partitioned into a clique $K$ and an independent set $S$. We suppose, without loss of generality, that $K$ is maximal, that is no vertex in $S$ is adjacent to all vertices in $K$. The pair $(K, S)$ is then called a canonical partition of $G$. For such a partition, we have $\omega(G) =|K|$.

We
follow \cite{MR1367739} for terminology and notation where they are not defined here.
The main tool we use in the proof of Theorem \ref{T4} is the Combinatorial
Nullstellensatz.

\begin{pro} (Combinatorial Nullstellensatz \cite{alon1992colorings})
Let $F$ be a field, let $d_1,\ldots, d_n \geq 0$ be integers, and let
$P \in F [x_1,\ldots, x_n]$ be a polynomial of degree $d_1 +\cdots +d_n$
with a non-zero $x_{1}^{d_{1}} \ldots x_{n}^{d_{n}}$ coefficient. Then $P$ cannot vanish on any set of the form $E_1 \times \cdots \times E_n$ with $E_1, \ldots, E_n \subset F$ and $|E_i| > d_i$ for $i = 1,\ldots, n$.
\end{pro}

\section{Proofs}

Here we prove that the difference between $dis[G]$ and $dis_{\ell}[G]$ can be arbitrary large.

\begin{ali}{
For every integer $t$, $t\geq 4$, we construct a graph $G$ such that $dis_{\ell}[G]-dis[G]\geq t$.
Our construction consists of four steps.\\ \\
{\bf Step 1.}
Consider $2t-1$ copies of the complete graph $K_{2t}$ and call them $\mathcal{K}_1,\mathcal{K}_2,\ldots,\mathcal{K}_{2t-1}$.
For every $i$, $1\leq i \leq 2t-1$, let $\{v_1^i,v_2^i,\ldots,v_t^i,u_1^i,u_2^i,\ldots,u_t^i\}$ be the set of vertices of the
complete graph $\mathcal{K}_i$. \\
{\bf Step 2.}
For each $(i,j,k)$, where $1\leq i < j \leq t$ and $ 1\leq k\leq 2t-1$, put two new vertices $x_{i,j}^k$ and $y_{i,j}^k$, and put the edges  $x_{i,j}^k y_{i,j}^k$, $x_{i,j}^k v_i^k$ and $ y_{i,j}^k v_j^k $.
Similarly, for every $(i,j,k)$, where $1\leq i < j \leq t$ and $ 1\leq k\leq 2t-1$, put two new vertices $a_{i,j}^k$ and $b_{i,j}^k$, and put the edges $a_{i,j}^k b_{i,j}^k$, $a_{i,j}^k u_i^k$ and $b_{i,j}^k u_j^k$. \\
{\bf Step 3.}
For every $(i,i',k)$, where $1\leq i \leq t$, $1\leq  i' \leq t$ and $ 1\leq k\leq 2t-1$, put two new  vertices $g_{i,i'}^k$ and $h_{i,i'}^k$, and put the edges $g_{i,i'}^k h_{i,i'}^k$, $g_{i,i'}^k v_i^k$ and
 $h_{i,i'}^k  u_{i'}^k$.\\
{\bf Step 4.}
Finally, put a new vertex $p$ and join the vertex $p$ to each vertex in $\{g_{i,i'}^k:  1\leq i \leq t,1\leq  i' \leq t,  1\leq k\leq 2t-1 \}$. Call the resulting graph $G$.
\\
\\
Next, we discuss the basic properties of the graph $G$. Let $f$ be a closed distinguishing labeling for $G$.
\\ \\
{\bf Lemma 1.1.} We have:\\
$d(v_{i}^k)=d(u_{i}^k)=\displaystyle 4t-2$,     \,\,\, for each $i,k$,   $1\leq i  \leq t$ and $ 1\leq k\leq 2t-1$,\\
$d(x_{i,j}^k)=d(y_{i,j}^k)=d(a_{i,j}^k)=d(b_{i,j}^k)=2$,\,\,\, for each $i,j,k$,   $1\leq i <j \leq t$ and $ 1\leq k\leq 2t-1$,\\
$d(g_{i,i'}^k)=3,\, d(h_{i,i'}^k)=2$,   \,\,\, for each $i,i',k$,   $1\leq i \leq t$, $1\leq  i' \leq t$ and $ 1\leq k\leq 2t-1$.
\\ \\
{\bf Lemma 1.2.} Let $M=\{ x_{i,j}^k,y_{i,j}^k,a_{i,j}^k,b_{i,j}^k:1\leq i <j \leq t, 1\leq k\leq 2t-1\}$. There is a function $f' : M\rightarrow \{1,2,\ldots,t\}$, such that for each $k$,
\begin{center}
$\displaystyle\sum_{l\in N[v_1^k]\cap M}f'(l), \ldots, \sum_{l\in N[v_t^k]\cap M}f'(l),\sum_{l\in N[u_1^k]\cap M}f'(l), \ldots, \sum_{l\in N[u_t^k]\cap M}f'(l)$
\end{center}
are $2t$ distinct integers.
\\ \\
{\bf Proof of Lemma 1.2.} Let $k$ be a fixed number and $f' : M\rightarrow \{1,2,\ldots,t\}$ be an arbitrary labeling.
For each $i$, $1\leq i \leq t  $ we have:
\begin{center}
$|\{l: l\in N[v_i^k]\cap M \} |=|\{l: l\in N[u_i^k]\cap M \} |= t-1$.
\end{center}
Thus,
\begin{center}
$t-1 \leq \displaystyle\sum_{l\in N[v_i^k]\cap M}f'(l) \leq t( t-1)$.
\end{center}
On the other hand, for each $i,j$, $i\neq j$, we have:
\begin{center}
$\{l: l\in N[v_i^k]\cap M \}\cap \{l: l\in N[v_j^k]\cap M \}=\emptyset$.
\end{center}
Also, for each $i,i'$, $1\leq i, i' \leq t$, we have:
\begin{center}
$\{l: l\in N[v_i^k]\cap M \}\cap \{l: l\in N[u_{i'}^k]\cap M \}=\emptyset$.
\end{center}
Since $t\geq 4$, one can define $f' : M\rightarrow \{1,2,\ldots,t\}$ such that for each $i$, $1\leq i \leq t  $, $$\sum_{l\in N[v_i^k]\cap M}f'(l)=t-1 +i-1$$ and $$\sum_{l\in N[u_{i}^k]\cap M}f'(l)=2t-1+i-1.$$ This completes the proof of Lemma $ \diamondsuit$
\\ \\
{\bf Lemma 1.3.} For each $(i,j,k)$, where $1\leq i < j \leq t$ and $ 1\leq k\leq 2t-1$,  $f(v_i^k)\neq f(v_j^k)$ and
$f(u_i^k)\neq f(u_j^k)$.
\\ \\
{\bf Proof of Lemma 1.3.}
Consider the two adjacent vertices $ x_{i,j}^k$ and $y_{i,j}^k$, since $f$ is a closed distinguishing labeling for $G$, we have,

\begin{center}
$\displaystyle\sum_{l \in N[x_{i,j}^k] }f(l) \neq \sum_{l \in N[y_{i,j}^k] }f(l)$.
\end{center}
Thus,
\begin{center}
$f(v_i^k) +f(x_{i,j}^k) + f(y_{i,j}^k) \neq f(v_j^k)+f(x_{i,j}^k) + f(y_{i,j}^k) $.
\end{center}
Therefore, $f(v_i^k)\neq f(v_j^k)$.
Similarly, by considering  the two adjacent vertices $ a_{i,j}^k$ and $b_{i,j}^k$, we have $f(u_i^k)\neq f(u_j^k)$. $ \diamondsuit$

By Lemma 1.3, $f(v_1^1), f(v_2^1), \ldots, f(v_t^1) $ are $t$ distinct integers. So $dis[G]\geq t$. Now, we show that
$dis[G]\leq t$. Let $f'$ be a labeling that has the conditions of Lemma 1.2 and consider the following labeling for $G$:
\\ \\
$f:V(G)\rightarrow \{1,2,\ldots,t\}$,
\\
$f(v_{i}^k)=f(u_{i}^k)=i$,     \,\,\, for each $i,k$,   $1\leq i  \leq t$ and $ 1\leq k\leq 2t-1$,\\
$f(g_{i,i'}^k)=f(h_{i,i'}^k)=1$,   \,\,\, for each $i,i',k$,   $1\leq i \leq t$, $1\leq  i' \leq t$ and $ 1\leq k\leq 2t-1$,\\
$f(p)=t$,\\
$f(l)=f'(l)$,   \,\,\, for each $l\in M$.\\ \\
Now, we show that $f$ is a closed distinguishing labeling for $G$. We have:\\\\
$\displaystyle\sum_{l \in N[p] }f(l)= t^2(2t-1)+t\geq 4t$, $\displaystyle\sum_{l \in N[g_{i,i'}^k] }f(l)= t+2+i  \leq 3t$, $\displaystyle\sum_{l \in N[h_{i,i'}^k] }f(l)= 2+i' \leq 3t$,\\ $\displaystyle\sum_{l \in N[x_{i,j}^k] }f(l)= f'(x_{i,j}^k)+f'(y_{i,j}^k)+i \leq 3t$,
$\displaystyle\sum_{l \in N[y_{i,j}^k] }f(l)= f'(x_{i,j}^k)+f'(y_{i,j}^k)+j \leq 3t$,\\
$\displaystyle\sum_{l \in N[a_{i,j}^k] }f(l)= f'(a_{i,j}^k)+f'(b_{i,j}^k)+i \leq 3t$,
$\displaystyle\sum_{l \in N[b_{i,j}^k] }f(l)= f'(a_{i,j}^k)+f'(b_{i,j}^k)+j \leq 3t$,\\
$\displaystyle\sum_{l \in N[v_{i }^k] }f(l)=  \displaystyle\sum_{l \in N[v_{i }^k] \cap M }f'(l)+ \displaystyle\sum_{l \in N[v_{i }^k] \setminus M }f(l) \geq 4t$,\\
$\displaystyle\sum_{l \in N[u_{i }^k] }f(l)=  \displaystyle\sum_{l \in N[u_{i }^k] \cap M }f'(l)+ \displaystyle\sum_{l \in N[u_{i }^k] \setminus M }f(l) \geq 4t$.\\ \\
For every two adjacent vertices $v_{i }^k$ and $u_{j }^k$, we have
$$\displaystyle\sum_{  l \in N[v_{i }^k] \setminus M  } f(l)=\displaystyle\sum_{  l \in N[u_{j }^k] \setminus M } f(l). $$
Thus, by Lemma 1.2,  the sum of labels of the vertices in the closed neighborhood of the vertex $v_{i }^k$ differs from
the sum of labels of the vertices in the closed neighborhood of the vertex $u_{j }^k$.
We have a similar result for every two adjacent vertices $v_{i }^k$ and $v_{j }^k$.
For other pairs of adjacent vertices, from the values shown above it is clear that for every two adjacent vertices $z$, $s$,
the sum of labels of the vertices in the closed neighborhood of the vertex $z$ differs from
the sum of labels of the vertices in the closed neighborhood of the vertex $s$.
So, $f$ is a closed distinguishing labeling for $G$.
Thus $dis[G]= t$.

Next, we show that $dis_{\ell}[G] \geq 2t$. To the contrary assume that $dis_{\ell}[G] \leq 2t-1$ and let $N=\{u_{i}^{k} :1\leq i  \leq t, 1\leq k\leq 2t-1\}$.
Consider the following lists for the vertices of the graph $G$:

$L(u_{i}^{k})=\{1+k,2+k,3+k, \ldots, 2t-1+k\}$,

$L(l)=\{1,2,3, \ldots, 2t-1\}$, for every $l\in V(G)\setminus N$.

Assume that $f$ is a closed distinguishing labeling for $G$ from the lists that shown above (i.e. for each vertex $v$, $f(v)\in L(v)$).
Without loss of generality assume that $f(p)=w$. Consider the set of vertices
$ v_1^w,v_2^w, \ldots,v_t^w, u_1^w,u_2^w, \ldots,u_t^w   $. We have:

$L(u_{i}^{w})=\{1+w,2+w,3+w, \ldots, 2t-1+w\}$,

$L(v_{i}^{w})=\{1,2,3, \ldots, 2t-1\}$.
\, \\ \\
Consider the following partition for the set of numbers $L(u_{i}^{w})\cup L(v_{i}^{w})$,

 \begin{center}
$\{1+w,1\},\{2+w,2\}, \ldots, \{2t-1+w,2t-1\}$.
 \end{center}

By the pigeonhole principle, there are indices $r$, $i$ and $j$ such that $f(v_{i}^{w}),f(u_{j}^{w})\in \{r+w,r\}$, so $f(v_{i}^{w})=r$ and $f(u_{j}^{w})=r+w$. Therefore,

\begin{center}
$\displaystyle\sum_{l\in N[g_{i,j}^{w}]}f(l)= \sum_{l\in N[h_{i,j}^{w}]}f(l) $.
 \end{center}
 This is a contradiction, so $dis_{\ell}[G] \geq 2t$. This completes the proof.

}\end{ali}

Here, we investigate the computational complexity of determining   $dis[G]$ for   planar subcubic graphs.
We show that for a given planar subcubic graph $G$, it is $ \mathbf{NP} $-complete to determine whether $ dis[G]=2$.

\begin{alii}{
 Let $\Phi$ be a $3$SAT formula with clauses $C=\lbrace
c_1, \ldots ,c_{\gamma}\rbrace $ and variables
$X=\lbrace x_1, \ldots ,x_n\rbrace $. Let $G(\Phi)$ be a graph with the vertices $C \cup X \cup (\neg X)$, where $\neg X = \lbrace \neg x_1, \ldots , \neg x_n\rbrace$, such that for each clause $c=(y \vee z \vee w )$, $c$ is adjacent to $y,z$ and $w$, also every $x \in X$ is adjacent to $\neg x$. $\Phi$ is called strongly planar formula if $G(\Phi)$ is a planar graph. It was shown that the problem of satisfiability for strongly planar formulas is $ \mathbf{NP}$-complete \cite{zhu1} (for more information about strongly planar formulas see \cite{dehghan2015strongly}). We reduce the following problem to our problem.

\textbf{Problem}: {\em Strongly planar $3$SAT.}\\
\textsc{Input}: A strongly planar formula  $ \Phi $.\\
\textsc{Question}: Is there a truth assignment for $ \Phi $ that satisfies all the clauses?\\

Consider an instance of strongly planar formula $\Phi$ with the variables
$X $ and the clauses $C $. We transform this into a planar subcubic graph $\mathcal{G}$
such that $ dis[\mathcal{G}]=2 $ if and only if $\Phi$ is satisfiable.
For every $x\in X$ consider a cycle $C_{24\gamma}$, where $\gamma$ is the number of clauses in $\Phi$ (call that cycle $\mathcal{C}_x$).
Suppose that $\mathcal{C}_x=v_1 v_2 \ldots  v_{24\gamma} v_1$ and color the vertices of $\mathcal{C}_x$ by function $\ell$.

$ \ell(v_i)=
\begin{cases}
   $red$,      &$     if  $\,\,  1\leq i \leq 12\gamma \, $ and $ i \equiv 1 \;(\bmod\; 6) \\
   $black$,    &$     if  $\,\,  1\leq i \leq 12\gamma   \, $ and $ i \equiv 4 \;(\bmod\; 6) \\
   $black$,    &$     if  $\,\,  1+12\gamma\leq i \leq 24\gamma \, $ and $ i \equiv 1 \;(\bmod\; 6) \\
   $blue$,     & $    if  $\,\, 1+12\gamma\leq i \leq 24\gamma  \,  $ and    $ i \equiv 4 \;(\bmod\; 6)\\
   $white$     & $    otherwise. $
\end{cases}$

For every $c\in C$
consider a path $P_{8}$ with the vertices $u_1,u_2,\ldots,u_8$, in that order. Next put two new isolated vertices $u'$ and $u''$, and join the vertex $u_3$ to the vertex $u'$ and join the vertex $u_6$ to the vertex $u''$. Call that resultant graph $\mathcal{P}_c$.
Next, for every $c\in C$, without loss of generality assume that $c=(a \vee b \vee w)$, where $a,b,w\in X \cup (\neg X)$.
If $a\in X$ ($a\in \neg X$) then join the vertex $u_1$, $u_1\in \mathcal{P}_c$ to one of the red (blue) vertices with degree two of $\mathcal{C}_a$. Similarly, if $b\in X$ ($b\in \neg X$) then join the vertex $u_1$, $u_1\in \mathcal{P}_c$ to one of the red (blue) vertices  with degree two of $\mathcal{C}_b$. Furthermore, if $w\in X$ ($w\in \neg X$) then join the vertex
$u_4$, $u_4 \in \mathcal{P}_c$ to one of the red (blue) vertices of degree two of $\mathcal{C}_w$; also, join the vertex
$u_8$, $u_8 \in \mathcal{P}_c$ to one of the red (blue) vertices of degree two of $\mathcal{C}_w$.
In the resulting graph for every red or blue vertex $l$ with degree two, put a new isolated vertex $l'$ and join the vertex $l$ to the vertex $l'$. Also, for every black vertex $l$, put a new isolated vertex $l'$ and join the vertex $l$ to the vertex $l'$. So in the final graph the degree of every blue, red or black vertex is three.
Call the resultant subcubic graph $\mathcal{G}$. Note that since $ \Phi $ is strongly planar ($G(\Phi)$ is planar),  we can construct $\mathcal{G}$
such that it is a planar graph.

Assume that $f:V(\mathcal{G})\rightarrow \{1,2\}$ is a closed distinguishing labeling for $\mathcal{G}$.
We have the following lemmas:
\\ \\
{\bf Lemma 2.1.} For every $x\in X$, we have:

$\triangleright$ for every $z\in V(\mathcal{C}_x) $, if $\ell(z)=$red then $f(z)=2$ and if $\ell(z)=$blue then $f(z)=1$,
\\
\\
or

$\triangleright$ for every $z\in V(\mathcal{C}_x) $, if $\ell(z)=$red then $f(z)=1$ and if $\ell(z)=$blue then $f(z)=2$.
 \\ \\
{\bf Proof of Lemma 2.1.} Let $ h_1, h_2,h_3,h_4, h_5,h_6, h_7\in V(\mathcal{C}_x)$, and $ h_1 h_2h_3h_4 h_5h_6 h_7$ be a path of length six in $\mathcal{C}_x$. Also, without loss of generality assume that
$\ell( h_2)=\ell( h_3)=\ell( h_5)=\ell( h_6)=$white. Since $f$ is a  a closed distinguishing labeling, we have:
\begin{center}
$\displaystyle\sum_{g\in N[h_2]}f(g)= \sum_{g\in N[h_3]}f(g) $.
 \end{center}
Thus, $f(h_1)\neq f(h_4)$. Similarly, $f(h_4)\neq f(h_7)$. Hence $ f(h_1)=f(h_7)$. Therefore, the labels of red vertices are the same. Also, the labels of blue vertices are the same. In $\mathcal{C}_x$, we have $\ell(v_1)=$red,  $\ell(v_{24\gamma -2})=$blue and $\ell(v_{24\gamma -1})=\ell(v_{24\gamma})=$white.  Thus, $f(v_{24\gamma -2})\neq f(v_1)$. This completes the proof. $ \diamondsuit$
\\
\\
Define $f': X \cup (\neg X) \rightarrow \{1,2\} $ such that  for every $a \in X \text{ }(a \in \neg X)$,  $f'(a)=2$ if and only if the values of function $f$ for the red (blue) vertices in $\mathcal{C}_a$ are two.
\\ \\
{\bf Lemma 2.2.} Let $c$ be an arbitrary clause and  $c=(a\vee b \vee w)$, where $a,b,w\in X \cup (\neg X)$.
We have $2\in \{f'(a),f'(b),f'(w)\}$.
\\ \\
{\bf Proof of Lemma 2.2.} To the contrary assume that $f'(a)=f'(b)=f'(w)=1$. Since
$\sum_{l\in N[u_1]}f(l)\neq \sum_{l\in N[u_2]}f(l)$,
we have $f(u_3)=1$. Also since $\sum_{l\in N[u_4]}f(l)\neq \sum_{l\in N[u_5]}f(l)$, we have $f(u_6)=1$. Finally, since
  $\sum_{l\in N[u_7]}f(l)\neq \sum_{l\in N[u_8]}f(l)$, we have $f'(w)=2$. But this is a contradiction.  $ \diamondsuit$

First, assume that $f:V(\mathcal{G})\rightarrow \{1,2\}$ is a closed distinguishing labeling for $\mathcal{G}$.
Let $\Gamma : X \rightarrow \{ true,false \} $ be a function such that  $\Gamma(x)=true$ if and only if $f'(x)=2$. By
Lemma 2.2, $\Gamma$ is a satisfying assignment for $ \Phi $.

Next, suppose that $ \Phi $ is satisfiable with the satisfying
assignment $ \Gamma $. For every $x\in X$ if $\Gamma(x)=true$ then for $ \mathcal{C}_x$ define:
$$ f(v_i)=
\begin{cases}
   2,      &$     if  $\,\,  \ell(v_i)=$red$ \\
   1,      &$     if  $\,\,  \ell(v_i)=$blue$ \\
   2,      &$     if  $\,\,  \ell(v_i)=$white$ \\
   1       &$     if  $\,\,  \ell(v_i)=$black$\,  $ and    $ 1 \leq i \leq 12\gamma \\
   2       &$     if  $\,\, \ell(v_i)=$black$   \,  $ and    $ 1+12\gamma \leq i \leq 24\gamma,
\end{cases}$$
and if $\Gamma(x)=false$ then for $ \mathcal{C}_x$ define:
$$ f(v_i)=
\begin{cases}
   1,       &$     if  $\,\,  \ell(v_i)=$red$ \\
   2,       &$     if  $\,\,  \ell(v_i)=$blue$ \\
   2,       &$     if  $\,\,  \ell(v_i)=$white$ \\
   2       &$     if  $\,\,  \ell(v_i)=$black$\,  $ and    $ 1 \leq i \leq 12\gamma \\
   1       &$     if  $\,\, \ell(v_i)=$black$   \,  $ and    $ 1+12\gamma \leq i \leq 24\gamma.
\end{cases}$$
Next, for every $c=(a \vee b \vee w)$, if $\Gamma(w)=true$ then for $ \mathcal{P}_c$, define:
$$ f(v_i)=
\begin{cases}
   2,       &$     if  $\,\,  v_i=u''  \\
   1       &$     otherwise  $,
\end{cases}$$
otherwise, if  $\Gamma(w)=false$ then for $ \mathcal{P}_c$, define:
$$f(v_i)=
\begin{cases}
   2,       &$     if  $\,\,  v_i\in\{u_3,u_6,u',u''\}  \\
   1       &$     otherwise  $.
\end{cases}$$
Finally, label  remaining vertices by number 2.
One can check that $f$ is   a closed distinguishing labeling for $\mathcal{G}$. This completes the proof.
}\end{alii}

Next, we show that it is $ \mathbf{NP} $-complete to determine  whether $ dis[G]=2$, for a given bipartite subcubic graph $G$.

\begin{aliii}{
We reduce {\em   Monotone Not-All-Equal 3Sat} to our problem in polynomial time. It was shown that the following problem is $ \mathbf{NP}$-complete \cite{MR1567289}.

 {\em Monotone Not-All-Equal 3Sat .}\\
\textsc{Instance}: Set $X$ of variables and collection $C$ of clauses over $X$ such that each
clause $c \in C$ has $\mid c  \mid = 3$ and there is no negation in the formula.\\
\textsc{Question}: Is there a truth assignment for $X$ such that each clause in $C$ has at
least one true literal and at least one false literal?\\

Consider an instance  $\Phi$ with the set of variables
$X $ and the set of clauses $C $. We transform this into a  bipartite graph $G$,  such that $\Phi$ has a Not-All-Equal satisfying assignment if and only if there is a closed distinguishing labeling $f:V(G) \rightarrow \{1,2\}$.
For every $x\in X$ consider a cycle $C_{12\gamma}$, where $\gamma$ is the number of clauses in $\Phi$ (call that cycle $\mathcal{C}_x$). Suppose that $\mathcal{C}_x=v_1v_2\ldots v_{12\gamma}v_1$ and color the vertices of $\mathcal{C}_x$ by function $\ell$.

$ \ell(v_i)=
\begin{cases}
   $red$,      &$     if  $ i \equiv 1 \;(\bmod\; 6) \\
   $blue$,    &$     if  $ i \equiv 4 \;(\bmod\; 6) \\
   $white$     & $    otherwise. $
\end{cases}$ \\ \\
For every $c=(x \vee y \vee z)$, $c\in C$, do the following three steps:\\ \\
{\bf Step 1.} Put two paths $\mathcal{P}_c=c_{1}^1 c_{2}^1 c_{3}^1 c_{4}^1 c_{5}^1 $ and $\mathcal{P}_c^{'}=c_{1}^2 c_{2}^2 c_{3}^2 c_{4}^2 c_{5}^2$. Also, put two isolated vertices $c'$, $c''$ and add the edges $c'c_{3}^1$, $c''c_{3}^2$.\\
{\bf Step 2.} Without loss of generality suppose that $\{ v_{i}, v_{j}, v_{k}\}$ is a set of vertices such that each of them has degree two, the value of function $\ell$ for each of them is red, $ v_{i}\in V(\mathcal{C}_x)$, $ v_{j}\in V(\mathcal{C}_y)$ and $v_{k}\in V(\mathcal{C}_z)$. Add the edges $ \displaystyle c_{1}^1 v_{i}, c_{1}^1 v_{j}, c_{5}^1 v_{k}$.\\
{\bf Step 3.} Without loss of generality suppose that $\{ v_{i'}, v_{j'}, v_{k'}\}$ is a set of vertices such that each of them has degree two, the value of function $\ell$ for each of them is blue, $ v_{i'}\in V(\mathcal{C}_x)$, $ v_{j'}\in V(\mathcal{C}_y)$ and $v_{k'}\in V(\mathcal{C}_z)$. Add the edges $ \displaystyle c_{1}^2 v_{i'}, c_{1}^2 v_{j'}, c_{5}^2 v_{k'}$.

Next, in the resulting graph for every red or blue vertex $u$ with degree two, put a new isolated vertex $u'$ and join the vertex $u$ to the vertex $u'$. Call the resultant bipartite subcubic graph $G$.

First, assume that $f$ is a closed distinguishing labeling for $G$.
For every $x\in X$, we have:

$\diamond$ for every $u\in V(\mathcal{C}_x) $, if $\ell(u)=$red then  $f(u)=2$, and if $\ell(u)=$blue then  $f(u)=1$,
\\ \\
or

$\diamond$ for every $u\in V(\mathcal{C}_x) $, if $\ell(u)=$red then  $f(u)=1$, and if $\ell(u)=$blue then  $f(u)=2$,
\\ \\
(see the proof of Lemma 2.1).
Define $\Gamma: X \rightarrow \{true,false\} $ such that  for every $x \in X$,  $\Gamma(x)=true$ if and only if the values of function $f$ for the red vertices in $\mathcal{C}_x$ are two. By the structure of clause gadgets, for
every clause $c=(x\vee y \vee z)$,

\begin{center}
$\displaystyle \sum_{u\in N[c_{1}^1]}f(u) \neq \sum_{u\in N[c_{2}^1]}f(u)$ and $\displaystyle \sum_{u\in N[c_{4}^1]}f(u) \neq \sum_{u\in N[c_{5}^1]}f(u)$.
 \end{center}
So, $true\in \{\Gamma(x),\Gamma(y),\Gamma(z)\}$.
On the other hand,

\begin{center}
$\displaystyle \sum_{u\in N[c_{1}^2]}f(u) \neq \sum_{u\in N[c_{2}^2]}f(u)$ and $\displaystyle \sum_{u\in N[c_{4}^2]}f(u) \neq \sum_{u\in N[c_{5}^2]}f(u)$.
 \end{center}
Thus, $false \in \{\Gamma(x),\Gamma(y),\Gamma(z)\}$.
Therefore, $\Gamma$ is a Not-All-Equal assignment.

Next, suppose that $ \Phi $ has a Not-All-Equal
assignment $ \Gamma $. For every $x\in X$ if $\Gamma(x)=true$ then:
\\ \\
$\diamond$ for every $u\in V(\mathcal{C}_x) $, if $\ell(u)=$red then put $f(u)=2$ and if $\ell(u)=$blue then put $f(u)=1$,
\\ \\
and if $\Gamma(x)=false$ then:
\\ \\
 $\diamond$ for every $u\in V(\mathcal{C}_x) $, if $\ell(u)=$red then put $f(u)=1$ and if $\ell(u)=$blue then put $f(u)=2$.
\\ \\
For every white vertex $l$, put $f(l) = 2$. Also, for every clause $c=(x \vee y \vee z)$, $c\in C$, put: $$f( c_{1}^1)=f( c_{2}^1)=f( c_{4}^1)=f( c_{5}^1)=f( c_{1}^2)=f( c_{2}^2)=f( c_{4}^2)=f( c_{5}^2)=1.$$ Also,
put $f( c_{3}^1) =1$ and $f( c_{3}^2) =2$ if and only if $\Gamma(z)=false$. Finally, label all remaining vertices by number 2. One can see that the resulting labeling is  a closed distinguishing labeling for $G$. This completes the proof.

}\end{aliii}

Here, we prove that for every integer $t\geq 3$, it is
{\bf NP}-complete to determine whether $dis[G]= t$ for a given graph
$G$.

\begin{aliiii}{
In order to prove the theorem, we reduce $t$-Colorability
to our problem for each $t\geq 3$. It was shown in \cite{MR1567289} that for each $t$, $t\geq 3$, the following problem is $ \mathbf{NP}$-complete.

\textbf{Problem}: {\em $t$-Colorability.}\\
\textsc{Input}: A graph $G$.\\
\textsc{Question}: Is $\chi(G)\leq t$?\\

Let $G$ be a given graph and $t$ be a fixed number. We construct a graph $G^*$ in polynomial time such that $\chi(G)\leq t$
if and only if $G^*$ has  a closed distinguishing labeling from $\{1,2,\ldots,t\}$.
Our construction consists of two steps.\\ \\
{\bf Step 1.}
Consider a copy of the graph $G$. For every vertex $v\in V(G)$ put $\Delta(G)-d_G(v)+1$ new isolated vertices $u_1^v,\ldots, u_{\Delta(G)-d_G(v)}^v,z^v$ and join them to the vertex $v$. Call the resulting graph $G'$. In the resulting graph the degree of each vertex is $\Delta(G)+1$ or $1$.\\
{\bf Step 2.}
Let $|V(G)|=n$, $|V(G')|=n'+n$ and $\alpha=(n'+1)(t-1)+2$. Consider a copy of the complete graph $K_{\alpha}$ with the set of  vertices $\{x_1,x_2,\ldots,x_{\alpha}\}$.
For each $i$, $1\leq i < \alpha $, put $n'+1$ new isolated vertices and join them to the vertex $x_i$. Finally,
for each $v\in V(G)$ join the vertex $x_{\alpha}$ to the vertices $v,u_1^v,\ldots, u_{\Delta(G)-d_G(v)}^v$.
Call the resulting graph $G^*$.
In the final graph for each $i$, $1\leq i < \alpha$, $d_{G^*}(x_i)=\alpha+n'$ and
$d_{G^*}(x_{\alpha})=\alpha+n'-1$.

Let $f:V(G^*)\rightarrow \{1,2,\ldots,t\}$ be a closed distinguishing labeling.
We have the following lemmas:
\\ \\
{\bf Lemma 4.1.} For every vertex $v\in V(G^*)$, $f(v)=f(u_1^v)=\cdots =f(u_{\Delta(G)-d_G(v)}^v)=1$.
\\ \\
{\bf Proof of Lemma 4.1.} Consider the set of vertices $V(K_{\alpha})=\{x_1,x_2,\ldots,x_{\alpha}\}$. Since $f$ is a a closed distinguishing labeling,

\begin{center}
$\displaystyle\sum_{l\in N[x_1]}f(l),  \sum_{l\in N[x_2]}f (l),  \ldots, \sum_{l\in N[x_{\alpha}]}f(l)$,
\end{center}
are $\alpha$ distinct numbers. Thus,

\begin{center}
$\displaystyle \sum_{l\in N[x_1] \atop l\notin V(K_{\alpha})}f(l),  \sum_{l\in N[x_2]\atop l\notin V(K_{\alpha}) }f (l),  \ldots, \sum_{l\in N[x_{\alpha}]\atop l\notin V(K_{\alpha})}f(l)$,
\end{center}
are $\alpha$ distinct numbers. For each $i$, $1\leq i < \alpha$,

\begin{center}
$\displaystyle n'+1 \leq \sum_{l\in N[x_i] \atop l\notin V(K_{\alpha})}f(l)\leq (n'+1)t$.
\end{center}
So,

\begin{center}
$\displaystyle\{\sum_{l\in N[x_i] \atop l\notin V(K_{\alpha})}f(l):1\leq i < \alpha \}=\{n'+1, n'+2, \ldots, (n'+1)t \}$.
\end{center}
Thus

\begin{center}
$\displaystyle \sum_{l\in N[x_{\alpha}] \atop l\notin V(K_{\alpha})}f(l)\leq n'$.
\end{center}
On the other hand,

\begin{center}
$\displaystyle \sum_{l\in N[x_{\alpha}] \atop l\notin V(K_{\alpha})}f(l)\geq |\{l:l\in N[x_{\alpha}], l\notin V(K_{\alpha})\}|= n'$.
 \end{center}
Therefore

\begin{center}
$\displaystyle \sum_{l\in N[x_{\alpha}] \atop l\notin V(K_{\alpha})}f(l)= n'$
 \end{center}
 and
for every vertex $v\in V(G^*)$, $f(v)=f(u_1^v)=\cdots =f(u_{\Delta(G)-d_G(v)}^v)=1$. This completes the proof of Lemma. $ \diamondsuit$
\\ \\
{\bf Lemma 4.2.} Let $v$ and $v'$ be two adjacent vertices in $G$. We have $f(z^v)\neq f(z^{v'})$.
\\ \\
{\bf Proof of Lemma 4.2.} For two adjacent vertices $v$ and $v'$ in $G$ we have $ \sum_{l\in N_{G^*}[v]}f(l)\neq
 \sum_{l\in N_{G^*}[v']}f(l)$. By Lemma 4.1 and since $d_{G^*}(v)=d_{G^*}(v')$, we have $f(z^v)\neq f(z^{v'})$. $ \diamondsuit$

Let $f:V(G^*)\rightarrow \{1,2,\ldots,t\}$ be a closed distinguishing labeling.
By Lemma 4.2, the following function is proper vertex $t$-coloring for $G$:
$$c:V(G)\rightarrow \{1,2,\ldots,t\}\text{ such that }c(v)=f(z^{v})$$
On the other hand, if $G$ is $t$-colorable (and $c$ is  a proper vertex $t$-coloring for the graph $G$), define:
\\ \\
 $f(v)=f(u_1^v)=\cdots =f(u_{\Delta(G)-d_G(v)}^v)=1$ and $f(z^{v})=c(v)$,  for every vertex $v\in V(G)$,\\
 $f(l)=t$,                                                                for every vertex $l\in T$ (note that $ T=\{x_1,x_2,\ldots,x_{\alpha}\}$).
\\ \\
Also, for every vertex $x_i$, $1\leq i < \alpha$, label the set of vertices $\{l: l\in N[x_{i}], l\notin V(K_{\alpha})\}$, such that in the final labeling
\begin{center}
$\displaystyle\{\sum_{l\in N[x_i] \atop l\notin V(K_{\alpha})}f(l):1\leq i < \alpha \}=\{n'+1, n'+2, \ldots, (n'+1)t \}$.
\end{center}
One can check that $f$ is a closed distinguishing labeling. This completes the proof.
}\end{aliiii}

In the next theorem we give some upper bounds by using the Combinatorial Nullstellensatz.

\begin{aliiiii}{
Let $V(G) = \{x_1, \ldots , x_n\}$ and $\displaystyle S_{i}:=\sum_{x_j\in N[x_i]}x_{j}$. Define the following polynomial:
$$\displaystyle f(x_{1},\ldots,x_n)=\prod_{x_t \sim x_s \atop t < s}^{N[x_t]\neq N[x_s]} (S_{t}-S_{s}).$$
\item[(i)]
Let $d=a_1+\cdots+a_n$ be the degree of $f(x_1,\ldots,x_n)$ and the coefficient of $x_{1}^{a_{1}}\ldots x_{n}^{a_{n}}$ be non-zero. Let $x_v$ be a vertex in $G$. The term $x_v$ appears in $S_t$ if $x_v\in N[x_t]$. Hence the term
$x_v$ appears in $(S_{t}-S_{s})$ if $x_v \in N[x_t] \cup N[x_s]$ and $x_v \notin N[x_t] \cap N[x_s]$.
Thus $x_v$ appears in $f(x_1,\ldots,x_n)$ at most $\sum_{x_j\in N(x_v)}(d(x_j)-1)$ times, which is less than or equal to  $(d_1-1)+\cdots+(d_{\Delta}-1)$. It means that the degree of $f(x_1,\ldots,x_{n})$ with respect to the variable $x_v$ is less than or equal to  $s$. Hence for each $i$, $1\leq i \leq n$, we have $a_{i}\leq s$.  Let $E_{1}=\cdots=E_n=\{1,2,\ldots,s+1\}$. By Combinatorial Nullstellensatz the polynomial $f(x_1,\ldots,x_{n})$ cannot vanish on the set $E_1 \times \cdots \times E_n$. It means that for each $i$, $1\leq i \leq n$, there exists $e_{i}\in E_{i}$ such that $f(e_1,\ldots,e_{n})\neq 0$. Then $e_1,\ldots,e_n$ is a closed distinguishing labeling for $G$.\\
\item[(ii)]
The degree of $f(x_1,\ldots,x_n)$ is $m$. Also, $f(x_1,\ldots,x_n)$ is not the zero polynomial. So for each monomial like $x_{1}^{a_1}\ldots x_{n}^{a_n}$ we have $a_{i}\leq m$. It is easy to check that there exists a monomial such that  $a_{i}< m$ for $1\leq i \leq n$. Therefore, the  Combinatorial Nullstellensatz finishes the proof.\\
\item[(iii)]
Since $\Delta -1\geq d_{t+1}\geq d_{t+2}\geq \cdots \geq d_{\Delta}$, it follows that $s\leq \Delta^{2}-2\Delta +t$. Also for $1\leq i \leq n$, $d_i\leq \Delta$. Hence $s\leq \Delta^2-\Delta$. \\
\item[(iv)]
Let $x_v$ be the only vertex such that $d(x_v)=\Delta$. Assign 1 as the label for the vertex $x_v$. For every $i$, the term  $x_i$ appears at most $(d_1-1)+\cdots+(d_{\Delta-1}-1)$ times. We have:
$$d_1+\cdots+d_{\Delta-1}-(\Delta-1)\leq \Delta +(\Delta-2)(\Delta-1)-(\Delta-1)=\Delta^2-3\Delta+3.$$
\item[(v)]
Let $x_v$ be a vertex in $G$. Let $N(x_v)=\{x_{a_{1}},\ldots,x_{a_{k}}\}$. The term $x_v$ doesn't appear in $(S_{a_{i}}-S_{b})$, where $x_v\in N(x_{a_{i}})\cap N(x_b)$. Then for every $1\leq i \leq n$, the term  $x_i$ appears at most $k(k-\lambda-1)$ times in $f(x_1,\ldots,x_n)$.  Then the  Combinatorial Nullstellensatz completes the proof.\\
\item[(vi)]
Let $x_v$ be a vertex in $G$ such that $d(x_v)=\alpha$. Assume that $N(x_v)=\{x_{a_{1}},\ldots,x_{a_{\alpha}}\}$. Let $V(G)-N[x_v]=\{x_{b_{1}},\ldots,x_{b_{n-\alpha-1}}\}$. Let $x_r,x_t\in V(G)$ be adjacent. Then $x_{v}$ appears in $S_r-S_t$ if and only if exactly one of $x_r,x_t$ belongs to $N(x_v)$ and another one belongs to $V(G)-N[x_v]$. So in $f$, the term $x_v$ appears at most $\alpha(n-\alpha-1)$, which is less than or equal to  $[\dfrac{n-1}{2}]^2$. Now the proof is complete by the Combinatorial Nullstellensatz.
}\end{aliiiii}

Here, we show that for each positive integer $t$ there is a bipartite graph $G$ such that $dis[G]>t$.

\begin{aliiiiii}{
Let $t$ be a fixed number. We construct a bipartite graph $G$ such that $dis[G]>t$.
Let $\alpha=t^2$.  Define:
\\ \\
$V(G)=X\cup Y \cup Z$,\\
$X=\{x_1,x_2,\ldots, x_{\alpha}\}$,\\
$Y=\{y_1,y_2,\ldots, y_{\alpha}\}$,\\
$Z=\{z_{A,B},z_{A,B}^{'}:A,B\subseteq \{1,2,3,\ldots,\alpha\}, A\neq \emptyset, B\neq \emptyset  \}$.
\\ \\
$\displaystyle E(G)=\bigcup_{A,B}\{z_{A,B}z_{A,B}^{'}, z_{A,B}x_i, z_{A,B}^{'}y_j: i\in A, j\in B\}$.

To the contrary assume that $f:V(G)\rightarrow \{1,2,3,\ldots,t\}$ is a closed distinguishing labeling.
For every two vertices $x_i$ and $y_j$, we have:
\begin{center}
$\displaystyle x_iz_{\{i\},\{j\}}, z_{\{i\},\{j\}}z'_{\{i\},\{j\}},y_jz'_{\{i\},\{j\}} \in E(G)$,
\end{center}
so
\begin{center}
$\displaystyle \sum_{l\in N[z_{\{i\},\{j\}}]}f(l) \neq \sum_{l\in N[z'_{\{i\},\{j\}}]}f(l)$,
\end{center}
thus $f(x_i)\neq f(y_j)$. Let $S_1$ and $S_2$ be two subsets of $\{1,2,3,\ldots,t\}$ such that $|S_1\cup S_2|\leq t$
and $S_1\cap S_2=\emptyset$. Without loss of generality, we can assume that for each $i$, $1\leq i \leq \alpha$,
$f(x_i)\in S_1$ and for every $j$, $1\leq j \leq \alpha$,
$f(y_j)\in S_2$.

Let $T_X=\{ f(x_1),\ldots,f(x_{\alpha}) \}$ and $T_Y=\{ f(y_1),\ldots,f(y_{\alpha}) \}$. We know that $f(x_i) \in S_1$ and $f(y_j) \in S_2$. Let $|S_{1}|=\mu$ and $|S_{2}|\leq t-\mu$. By the pigeonhole principle there exists an element, say $r$, in $S_1$ such that $r$ appears at least $\dfrac{\alpha}{\mu}$ times in $T_X$. Similarly, there exists an element, say $p$, in $S_2$ such that $p$ appears at least $\dfrac{\alpha}{t-\mu}$ times in $T_Y$. Since $\alpha=t^2$, we have:
$$\dfrac{\alpha}{\mu}= \dfrac{t^2}{\mu}\geq \dfrac{t^2}{t}\geq p .$$
Also,
$$\dfrac{\alpha}{t-\mu}= \dfrac{t^2}{t-\mu}\geq \dfrac{t^2}{t}\geq r .$$
Thus in $T_X$ there exists $r$ at least $p$ times, and in $T_Y$ there exists $p$ at least $r$ times.
Consequently, one can find two sets $A,B\subseteq \{1,2,3,\ldots,\alpha\}$ such that $|A|=p$, $|B|=r$, for each $i\in A$, $f(x_i)=r$ and for each $j\in B$, $f(y_j)=p$. Thus,
\begin{center}
$\displaystyle \sum_{l\in N[z_{A,B}]}f(l) =\sum_{l\in N[z'_{A,B}]}f(l)$.
\end{center}
But this is a contradiction. Therefore, $dis[G]>t$.

}\end{aliiiiii}

Note that in the graph $G$ which was constructed in the previous theorem, we have:
$$|V(G)|=2\alpha+2(2^{\alpha}-1)^2=2t^2 + 2(2^{t^2}-1)^2.$$
It is interesting to find a bipartite graph $G$ such that $V(G)= \mathcal{O}(t^c)$ and $dis[G]>t$, where $c$ is a constant number.
Next, we show that if $G$ is a split graph, then $dis[G]\leq(\omega(G))^2$.

\begin{aliiiiiii}{
Let $(K,S)$ be a canonical partition of $G$ and assume that $S=\{v_1,v_2,\ldots,v_{|S|}\}$.
For each $i$, $1 \leq i \leq |S|$, let  $G_i$ be the induced subgraph on the set of vertices $K \bigcup_{j=1}^{i} v_j$.
Let  $G_0$ be the induced graph on the set of vertices $K$ and $f_0:V(G_0)\rightarrow \{1\}$ be a closed distinguishing labeling such that for every vertex $u\in V(G_0)$, $f_0(u)=1$. For $i=1$ to $i=|S|$ do the following procedure:\\ \\
$\triangleright$ For each $j$, $1 \leq j  \leq (\omega(G))^2$, define $$ g_{j}^{i}(x)=
     \begin{cases}
       f_{i-1}(x) & x\in V(G_{i-1}), \\
       j          & x\in V(G_i)\setminus V(G_{i-1}).\\
     \end{cases}$$
For a fixed number $j$ if there are two vertices $x,y\in V(G_0)$, such that $N_{G_i}(x)\neq N_{G_i}(y)$ and $$\displaystyle\sum_{l\in N_{G_i}[x]}g_j^{i}(l)=\sum_{l\in N_{G_i}[y]}g_j^{i}(l), $$ then   $g_{j}^{i}(x)$ is not a closed distinguishing labeling for
$G_i$.
The graph $G_0$ has $\displaystyle {{\omega(G)}\choose {2}}$ edges, so there are at most $\displaystyle {{\omega(G)}\choose {2}}$ restrictions. Thus, there is an index $j$ such that for every two vertices $x,y\in V(G_0)$, if $N_{G_i}(x)\neq N_{G_i}(y)$, then $$\displaystyle\sum_{l\in N_{G_i}[x]}g_j^{i}(l)\neq\sum_{l\in N_{G_i}[y]}g_j^{i}(l). $$
For that $j$, put $f_i \leftarrow g_j^{i} $. (End of procedure.)

When the  procedure terminates, the function $f_{|S|}$ is a closed distinguishing labeling for $G$. This completes the proof.
 }\end{aliiiiiii}

Here, we show that for each $n$, there is a graph  $G$ with $n$ vertices such that $dis[G]-dis_{s}[G]=\Omega(n^{\frac{1}{3}})$.

\begin{alij}{
Let $t=10k$ and consider a copy of the complete graph $K_{t^2}$ with the set of vertices $\{v_{i,j}: 1\leq i \leq t, 1\leq j \leq t,\}$.
For each $(i,j)$, $1\leq i \leq t, 1\leq j \leq t$, put $i+j$ new vertices $x_1^{i,j},x_2^{i,j},\ldots,x_i^{i,j},y_1^{i,j},y_2^{i,j},\ldots,y_j^{i,j}$ and join them to the vertex $v_{i,j}$.
Call the resultant graph $G$. Note that for each $(i,j)$, $1\leq i \leq t, 1\leq j \leq t$, $d(v_{i,j})=t^2+i+j-1$.
\\ \\
First, we show that $dis_{s}[G]=2$. Define:

$f:V(G)\rightarrow\{1,\Delta(G)+1\}$,

$f(v_{i,j})=f(x_1^{i,j})=f(x_2^{i,j})=\cdots=f(x_i^{i,j})=\Delta(G)+1$, for every $i$ and $j$,

$f(y_1^{i,j})=f(y_2^{i,j})=\cdots=f(y_j^{i,j})=1$, for every $i$ and $j$.
\\ \\
It is easy to check that $f$ is  a closed distinguishing labeling for $G$. Next, assume that
$f:V(G)\rightarrow dis[G]$ is  a closed distinguishing labeling for $G$.
Consider the set of vertices $R=\{v_{i,j}: 1\leq i \leq t, 1\leq j \leq t\}$. The function $f$ is  a closed distinguishing labeling therefore,

\begin{center}
$\displaystyle\{\sum_{l\in N[v_{i,j}]}f(l):  1\leq i \leq t, 1\leq j \leq t\}$,
\end{center}
are $t^2$ distinct numbers. Thus,

\begin{center}
$\displaystyle\{\sum_{l\in N[v_{i,j}] \atop l\notin R}f(l):  1\leq i \leq t, 1\leq j \leq t\}$,
\end{center}
are $t^2$ distinct numbers. For each $(i,j)$, $1\leq i \leq t, 1\leq j \leq t$,

\begin{center}
$\displaystyle 2\leq i+j \leq \sum_{l\in N[v_{i,j}] \atop l\notin R}f(l) \leq (i+j)dis[G] \leq (2t)\times dis[G]$.
\end{center}
So $2t\times dis[G]-2+1 \geq t^2$. Thus $dis[G]\geq 5k$.
On the other hand, $$\displaystyle|V(G)|=t^2+\sum_{i=1}^{t}\sum_{j=1}^{t}(i+j)\leq t^2+\sum_{i=1}^{t}\sum_{j=1}^{t}(2t)=\mathcal{O}(t^3)=\mathcal{O}(k^3).$$
This completes the proof.
}\end{alij}

Here, we show that for a given 4-regular graph $G$, it is $ \mathbf{NP} $-complete to determine  whether $ dis_{s}[G]=2$.

\begin{alijj}{
Clearly, the problem is in $ \mathbf{NP} $. We prove the $ \mathbf{NP} $-hardness by a reduction from the  following well-known $ \mathbf{NP} $-complete problem \cite{MR1567289}.

 {\em 3SAT.}\\
\textsc{Instance}: A 3CNF formula $\Psi=(X,C)$.\\
\textsc{Question}: Is there a truth assignment for $X$?

\begin{figure}[ht]
\begin{center}
\includegraphics[scale=.4]{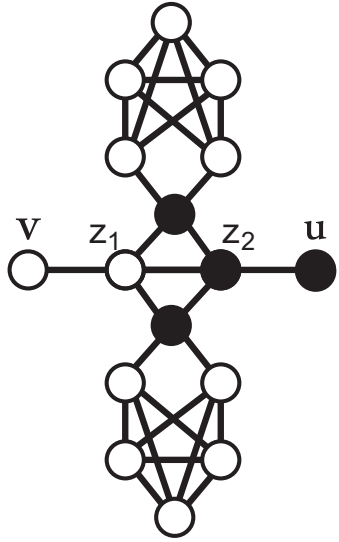}
\caption{The gadget $I(v,u)$. Let $G$ be a 4-regular graph and $f:V(G)\rightarrow \{\alpha,\beta\}$ be a closed distinguishing labeling.
If $I(v,u)$ is a subgraph of $G$, then $f(v)\neq f(u)$.
} \label{P1}
\end{center}
\end{figure}

Let $\Psi=(X,C)$ be an instance of {\em 3SAT} and also assume that $\alpha $ and $\beta$ are two numbers such that $\alpha\neq \beta$. We convert $\Psi$ into a 4-regular graph $G$ such that
$\Psi$ has a satisfying assignment if and only if $G$ has a closed distinguishing labeling from $\{\alpha ,\beta\}$.
First, we introduce a useful gadget.

{\bf Construction of the  gadget $\mathcal{T}_k$.}\\
Consider a copy of the bipartite graph $P_{2} \square C_{2k}$ and let $\ell :V(P_{2} \square C_{2k})\rightarrow \{1,2\}$ be a proper vertex 2-coloring. Call the set of vertices $V(P_{2} \square C_{2k})$, {\it the main vertices}. Construct the gadget $\mathcal{T}_k$ by replacing every edge $vu$ of $P_{2} \square C_{2k}$ with a copy of the gadget $I(v,u)$ which is shown in Fig. 1.

Note that the gadget $\mathcal{T}_k$ has $ 4k$ main vertices and the degree of each main vertex is three. Also, in
$\mathcal{T}_k$ the degree of each vertex that is not a main vertex is four.

For each variable $x\in X$ assume that the variable $x$ appears in exactly $\mu(x)$ clauses (positive or negative)
and suppose that  $|C|=\lambda$.  Next, we present the construction of the main graph.

{\bf Construction of the graph $G$.}\\
Put a copy of $\mathcal{T}_{3\lambda}$ and call it $F$. Also, for every variable $x\in X$, put a copy of the gadget $\mathcal{T}_{\mu(x)}$ and call it $D_x$. Furthermore, for every
clause $c\in C$, put a copy of the path $P_2=c_1c_2$.
For every $x\in X$, define:

$$ S_x^2 = \{  v\in V(D_x):  v\text{ is a main vertex and }  \ell(v)=2\}, $$

$$ S_x^1 = \{  v\in V(D_x):  v\text{ is a main vertex and } \ell(v)=1\}. $$
Also, define:

$$ R^2 = \{  v\in V(F):  v\text{ is a main vertex and } \ell(v)=2\}, $$

$$ R^1 = \{  v\in V(F):  v\text{ is a main vertex and } \ell(v)=1\}. $$
Next, for every $c\in C$, without loss of generality assume that $c=(a\vee b\vee s)$, where $a,b,s\in X \cup (\neg X)$.
If $a\in X$ ($a\in \neg X$) then join the vertex $c_1$, to a vertex $v\in S_a^2$ ($v\in S_{\neg a}^1$) of degree three.
Also, if $b\in X$ ($b\in \neg X$) then join the vertex $c_1$, to a vertex $v\in S_b^2$ ($v\in S_{\neg b}^1$) of degree three.
Similarly, if $s\in X$ ($s\in \neg X$) then join the vertex $c_1$, to a vertex $v\in S_s^2$ ($v\in S_{\neg s}^1$) of degree three.\\
Furthermore, join the vertex $c_2$ to three vertices
$v,u,z\in R^1$ of degree three.
Call the resultant graph $G'$. Note that the degree of every vertex in $G'$ is three or four.

Now, consider two copies of the graph $G'$. For each vertex $ h$ with degree three in $G'$, call its corresponding vertex in the first copy of $G'$, $h'$, and call its corresponding vertex in the second copy of $G'$, $h''$. Next, connect the vertices $h'$ and $h''$ through a copy of the gadget $I(h',h'')$.
Call the resulting 4-regular graph $G$. In the next, we just focus on the vertices in the first copy of $G'$ and talk about their properties.

First, assume that $f:V(G)\rightarrow \{\alpha,\beta\}$ is a closed distinguishing labeling. We have the following lemmas about the vertices in the first copy of $G'$.
\\ \\
{\bf Lemma 9.1.} For each $x\in X$, for every two vertices $h,g \in S_x^2$, $f(h)=f(g)$
and for every two vertices $h,g \in S_x^1$, $f(h)=f(g)$. Also, for each two vertices $h\in S_x^2$ and $g\in S_x^1$, $f(h)\neq f(g)$.
\\ \\
{\bf Proof of Lemma 9.1.} Let $G$ be a 4-regular graph and $f:V(G)\rightarrow \{\alpha,\beta\}$ be a closed distinguishing labeling for $G$. Assume that  $I(v,u)$ is a subgraph of $G$. For two adjacent vertices $z_1  $ and $z_2$ in $I(v,u)$ we have:
\begin{center}
$\displaystyle\sum_{l\in N[z_1]}f(l) \neq \displaystyle\sum_{l\in N[z_2]}f(l)$.
\end{center}
Thus, $f(v)\neq f(u)$. Consequently, in each copy of the gadget $I(v,u)$, we have $f(v)\neq f(u)$.
So, in the gadget $D_x$, for every two main vertices $l_1  $ and $l_2$ that are connected through a
copy of gadget  $I(l_1,l_2)$, we have $f(l_1)\neq f(l_2)$. On the other hand, the gadget $D_x$ is constructed from a bipartite graph by replacing each edge with a copy of the gadget $I(v,u)$. The main vertices of
$D_x$ can be partitioned into two sets, based on the function   $\ell$ which is  a proper vertex 2-coloring for the base bipartite graph.
In each part, the values of function $f$ for the main vertices in that part are the same.
So, for each $x\in X$, for every two vertices $h,g \in S_x^2$, $f(h)=f(g)$
and for every two vertices $h,g \in S_x^1$, $f(h)=f(g)$. Also, for each two vertices $h\in S_x^2$ and $g\in S_x^1$, $f(h)\neq f(g)$. $ \diamondsuit$
\\ \\
{\bf Lemma 9.2.} For every two vertices $g,h \in R^2$, $f(g)=f(h)$ and for every two vertices $g,h \in R^1$, $f(g)=f(h)$.
Also, for each two vertices $v\in R^2$ and $u\in R^1$, $f(g)\neq f(h)$.
\\ \\
{\bf Proof of Lemma 9.2.} The proof is similar to the proof of Lemma 9.1. $ \diamondsuit$
\\ \\
Note that if $f:V(G)\rightarrow \{\alpha,\beta\}$ is a closed distinguishing labeling, then
$$f'(v)=\begin{cases}
   \beta,      &$     if  $\,\,  f(v)=\alpha, \\
   \alpha,      &$     if  $\,\,  f(v)=\beta.
\end{cases}$$
is a closed distinguishing labeling for $G$.
Now, without loss of generality assume that the values of function $f$ for the set of vertices $R^2$ are $\beta$.
Define $\Gamma: X \rightarrow \{true,false\} $ such that for every $x \in X$,  $\Gamma(x)=true$ if and only if the values of function $f$ for the set of vertices $S_x^2$ are $\beta$.
By Lemma 9.1, Lemma 9.2 and structure of $G$,
it is easy to check that $\Gamma$ is a satisfying assignment. \\ \\
Next, suppose that $ \Psi $ is satisfiable with the satisfying assignment $ \Gamma $. Define the following values for the function $f$ for the vertices in the first copy of $G'$:
\\ \\
$\triangleright$ For every $v\in  S_x^2 $, if $\Gamma(x)=true $ then put $f(v)=\beta$ and if $\Gamma(x)=false $ then put $f(v)=\alpha$.
\\ \\
$\triangleright$ For every $v\in  S_x^1 $, if $\Gamma(x)=true $ then put $f(v)=\alpha$ and if $\Gamma(x)=false $ then put $f(v)=\beta$.
\\ \\
$\triangleright$ For every $v\in  R^2 $,  put $f(v)=\beta$ and for every $v\in  R^1 $,  put $f(v)=\alpha$.
\\ \\
$\triangleright$ For every $c\in  C $,  put $f(c_1)=\beta$ and $f(c_2)=\alpha$.
\\ \\
Also, for every vertex $l$, $l\in S_x^2 \cup S_x^1 \cup R^2 \cup R^1 \cup \{c_1,c_2: c\in C\}$ in the second copy of $G'$, put
$f(l)=1$ if and only if the value of function $f$ for the vertex $l$ in the first copy of $G'$ is two.
Finally, for each subgraph $I(v,u)$, without loss of generality assume that $f(v)=\alpha$ and $f(u)=\beta$. Label the vertices of
$V(I(v,u))\setminus\{v,u\}$ such that the labels of black vertices are $\beta$ and the labels of white vertices are $\alpha$ (see Fig. 1).
It is easy to check that this labeling is   a closed distinguishing labeling for $G$. This completes the proof.
}\end{alijj}

\section{Concluding remarks and future work}

In this article, we worked on
the closed distinguishing labeling which is very similar to the concept of relaxed locally identifying coloring.
A  vertex-coloring of a graph $G$ (not necessary proper) is said to be  relaxed locally identifying if for any pair $u$, $v$ of adjacent vertices with distinct closed neighborhoods, the sets of colors in the closed neighborhoods of $u$ and $v$ are different and an assignment of numbers to the
vertices of graph $G$ is {\it closed distinguishing}
if for any two adjacent vertices $v$ and $u$ the sum of labels of the vertices in the closed neighborhood of the vertex $v$ differs from
the sum of labels of the vertices in the closed neighborhood of the vertex $u$
unless  they have the same closed neighborhood.

\subsection{The computational complexity}

We proved that for a given bipartite subcubic graph $G$, it is $ \mathbf{NP} $-complete to decide  whether $ dis[G]=2$.
On the other hand, it was shown that for every tree $T$, $dis_{\ell}[T]\leq 3$ \cite{dis}. Here, we ask the following question.

\begin{prob}
For a given tree $T$, for every vertex $v\in V(T)$, let $L(v)$ be a list of size two of natural numbers.
Determine the computational complexity of deciding whether there is a closed distinguishing labeling $f$ such that for each $v \in V(T)$, $f(v) \in L(v)$.
\end{prob}

It was shown in  \cite{dis} that for every tree $T$, $dis[T]\leq 2$.
On the other hand, we proved that for a given bipartite  subcubic graph $G$ it is $ \mathbf{NP} $-complete to decide whether $ dis[G]=2$. In the proof of Theorem \ref{T3}, we reduced Not-All-Equal to our problem and the planar version of Not-All-Equal
is in $ \mathbf{P} $ \cite{p}, so the computational complexity of
deciding whether $ dis[G]=2$ for planar bipartite graphs remains unsolved.

\begin{prob} \label{pb}
For a given planar bipartite graph $G$, determine the computational complexity of deciding whether $dis[G]=2$.
\end{prob}

Let $G$ be an $r$-regular graph. If $dis[G]=2$ then for every two numbers $a,b$ $(a\neq b)$, $G$ has a closed distinguishing labeling from $\{a,b\}$.
We proved that for a given 4-regular graph $G$, it is $ \mathbf{NP} $-complete to decide  whether $ dis[G]=2$.
Determining the computational complexity of
deciding whether $ dis[G]=2$ for 3-regular graphs can be interesting.

\begin{prob} \label{3reg}
For a given 3-regular graph $G$, determine the computational complexity of deciding whether $dis[G]=2$.
\end{prob}

Summary of results and open problems in the complexity of determining whether $dis[G]=2$ is shown in Table 1.

\begin{table}[ht]
\small
\caption{Recent results and open problems} 
\centering 
\begin{tabular}{| l | c | c | } 
\hline\hline 
                 & \, \,\, $dis[G]=2$? \, \,\, & \, \,\, Refrence \, \, \, \\ \hline
\hline 
Tree &  {\bf P} &  see \cite{dis}  \\  \hline
Planar bipartite &  Open  &  Problem \ref{pb}  \\ \hline
Bipartite subcubic \, \, &    {\bf NP}-c &  Theorem \ref{T3}  \\ \hline
Planar subcubic  &     {\bf NP}-c &  Theorem  \ref{T2} \\ \hline
3-regular  &  Open &  Problem \ref{3reg}  \\ \hline
4-regular  &  {\bf NP}-c &  Theorem  \ref{T7}   \\
 [1ex] 
\hline 
\end{tabular}
\label{table:nonlin} 
\end{table}

\subsection{Bipartite graphs}

Let $G$ be a bipartite graph with partite sets $A$ and $B$ which is not a star. Let,
for $X\in \{A,B \}$; $\displaystyle \Delta_{X} = \max_{x\in X}d(x) $
and $\displaystyle\delta_{X,2} = \min_{x\in X, d(x)\geq 2}d(x)$. It was shown in \cite{dis} that
$$\displaystyle dis[G] \leq min\{c\sqrt{|E(G)|},\lfloor \dfrac{\Delta_A -1}{\delta_{B,2}-1}\rfloor+1,  \lfloor \dfrac{\Delta_B -1}{\delta_{A,2}-1}\rfloor+1\}, $$
where $c$ is some constant. Thus, for a given bipartite graph $G$, $dis[G]=\mathcal{O}(\Delta)$ \cite{dis}.
On the other hand, we proved that for each integer $t$, there is a bipartite graph $G$ such that $dis[G]>t$ (to see an example see Fig. 2). Here, we ask the following:

\begin{prob}
For each positive integer $t$, is there a bipartite graph $G$ such that $V(G)= \mathcal{O}(t^c)$ and $dis[G]>t$, where $c$ is a constant number.
\end{prob}

What can we say about the upper bound in bipartite graphs?
Perhaps
one of the most intriguing open question in this scope is the case of bipartite graphs.

\begin{prob}
Let  $G$ be a bipartite graph, is $dis[G]\leq \mathcal{O}(\sqrt{\Delta(G)})$?
\end{prob}

\begin{figure}[ht]
\begin{center}
\includegraphics[scale=.4]{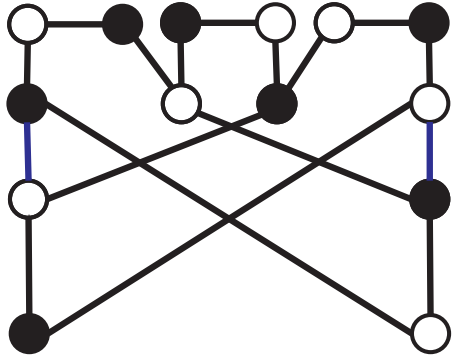}
\caption{A bipartite graph $G$ with the closed distinguishing number greater than two.
} \label{P1}
\end{center}
\end{figure}

We proved that the difference between $dis[G]$ and $dis_{\ell}[G]$ can be arbitrary large.
What can we say about the difference in bipartite graphs?

\begin{prob}
For any positive integer $t$, is there any bipartite graph $G$ such that $ dis_{\ell}[G]-dis[G] \geq t$?
\end{prob}

\subsection{General upper bounds and lower bounds}

For a given bipartite graph $G=[X,Y]$, define $f:V(G)\rightarrow \{1,\Delta\}$ such that:
$$ f(v)=
\begin{cases}
   1,           &$     if  $ v\in X \\
   \Delta,    &$     if  $ v\in Y
\end{cases}$$
It is easy to see that $f$ is a closed distinguishing labeling for $G$. Thus, for a bipartite graph $G$, $dis_s[G]\leq 2 $. On the other hand, for a general graph $G$, the best upper bound we know is
$dis_s[G]\leq |V(G)|  $.

\begin{prob}
Is this true "for any graph $G$,  $dis_s[G]\leq \chi(G)$?"
\end{prob}

For each positive integer $n$, we proved that there is a graph  $G$ with $n$ vertices such that $dis[G]-dis_{s}[G]=\Omega(n^{\frac{1}{3}})$. It would be desirable to increase the gap into $\Omega(\sqrt{n})$.

\begin{prob}
Is this true? "For each positive integer $n$, there is a graph  $G$ with $n$ vertices such that $dis[G]-dis_{s}[G]=\Omega(\sqrt{n})$".
\end{prob}

\subsection{Split graphs}

It is well-known that split graphs can be recognized in polynomial time, and that finding a canonical partition of a split graph can also be found in polynomial time. In this work, we proved that if $G$ is a split graph, then $dis[G]\leq(\omega(G))^2$.
Let $G$ be a split graph and $(K,S)$ be a canonical partition of $G$. Assume that $S=\{v_1,v_2,\ldots,v_{|S|}\}$.
Define:
$$ f(u)=
\begin{cases}
   1,                          &$     if  $\,\,  u\in V(K) \\
   (\Delta+1)^{i-1},           &$     if  $\,\,  u=v_i \text{ and } 1\leq i \leq |S|
\end{cases}
$$
It is easy to check that $f$ is a closed distinguishing labeling for $G$. Thus, $dis_s[G]\leq \alpha(G)$.
However, one further step does not seem trivial.

\begin{prob}
Is it true that if $G$ is a split graph, then $dis[G]=\mathcal{O}(\omega(G))$?
\end{prob}

\begin{prob}
Can one decide in polynomial time whether $dis[G]\leq \omega(G)$ for every split graph $G$?
\end{prob}

\section{Addendum}

During the review of the paper, Axenovich et al. in \cite{dis} independently, proved Theorem \ref{Tj1} and put it in the final version of their work.

\section{Acknowledgment}

The authors would like to thank the anonymous referees for their useful comments
and suggestions, which helped to improve the presentation of this paper.
The  work of the first author was done during a sabbatical at the School of Mathematics and Statistics, Carleton University, Ottawa. The first author is grateful to   Brett Stevens for hosting this
sabbatical.

\bibliographystyle{plain}
\bibliography{DISref}

\end{document}